\documentclass{article}
\usepackage{amsfonts}
\usepackage{bbm}
\usepackage{mathrsfs}
\usepackage{amsfonts}
 \input amssymb.sty
 \usepackage[notref,notcite]{showkeys}
 \usepackage{mathrsfs,amsfonts,amsmath,amssymb}
\usepackage[dvips]{color}
 \usepackage{amssymb}
 \usepackage[dvips]{color}
\allowdisplaybreaks
 \newtheorem{lemma}{\noindent\mbox{Lemma}}[section]
 \newtheorem{theorem}{\noindent\mbox{Theorem}}[section]

 \def\qed{\hfill$\Box$\medskip}
 \def\rto{\rightarrow\infty}
 \def\z{\left}
 \def\y{\right}
 \def\no{\nonumber}

\begin{document}
 \title{{\bf On the slowdown  of random
walk in random environment with bounded jumps}\thanks{The project is
partially supported by National Nature Science Foundation of
China(Grant No.11226199).}}

\author{ Huaming Wang\thanks{Department of Mathematics, Anhui Normal University, Wuhu 241000, P. R. China. hmking@mail.ahnu.edu.cn} }
\date{}
\maketitle \quad\\
\vspace{-6mm}

\begin{center}
\begin{minipage}[c]{14cm}
\begin{center}\textbf{Abstract}\quad \end{center}
\vspace{-0em} \hspace{2em} In this paper we prove that under certain
assumptions the transient random walk in random environment with
bounded jumps (in $\mathbb{Z}$) grows much slower than the speed
$n$. Precisely, there is $0<s<1$, such that although $X_n\rto$ we
have
 $\frac{X_n}{n^{s'}}\rightarrow 0$  for $0<s<s'$ almost surely. \\
\mbox{}\hspace{2.3em}\textbf{Keywords:}\quad random walk in random environment, slowdown, trap.\\
\mbox{}\hspace{2.3em}\textbf{AMS Subject Classification: \quad
Primary 60F15; secondary 60F10}.
\end{minipage}
\end{center}

\begin{center}
{\section{\hspace{-0.5cm.}  Introduction}}
\end{center}
Slowdown property is  one of the most important feature for the
random walk in random environment (RWRE in short) in $\mathbb{Z}.$
More precisely, although $X_n\rto$ we have that
$\frac{X_n}{n}\rightarrow 0$ almost surely (\cite{Z04}).  However,
this phenomena is impossible for random walk in non-random
environment, since the law of large numbers implies that the walk
grows with a positive speed as long as it is transient. Intuitively,
because of the random environment, there are ``many"  environments
formulated ``traps" in which the random walk spend ``much" time. For
the nearest RWRE (i.e., the walk which goes to right and left for
only one unit in one step), even a much slower speed has been
revealed, i.e., under certain assumptions there is $0<s<1$, such
that although $X_n\rto$ we have
 $\frac{X_n}{n^{s'}}\rightarrow 0$  for $0<s<s'$ almost surely (\cite{S04}).
In this paper, we will prove this  property  for the random walk in
random environment with bounded jumps. One should note that the RWRE
with bounded jumps makes the situation more complicated than the
nearest RWRE (\cite{B02}).

Let us recall the RWRE with bounded jumps firstly. We will adapt the
notations in \cite{B02}. $\Lambda=\{-L,...,1\} $,  $\Sigma$ is the
simplex in $ \mathbb{R}^{L+2},$ and $\Omega:=\Sigma^\mathbb{Z}.$ Let
$\mu$ be a measure on $\Sigma$ and
$\omega_0=(\omega_0(z))_{z\in\Lambda}$ be a $\Sigma$-valued random
vector with distribution $\mu,$ satisfying
$\sum_{z\in\Lambda}\omega_0(z)=1,$ and
$\mu(\omega_0(z)/\omega_0(1)>\kappa,\ z\in \Lambda,\ z\ne 0)=1 $ for
some $\kappa>0.$ Let $\mathbb{P}=\mu^{\bigotimes\mathbb{Z}}$ on
$\Omega$ making $\omega_x, x\in\mathbb{Z}$ i.i.d.. The random walk
in random environment $\omega$ with bounded jumps is the Markov
chain defined by $X_0=x$ and the transition probabilities
\begin{equation}\label{m1}
    P_{x,\omega}(X_{n+1}=y+z|X_n=y)=\omega_y(z), \forall y\in
    \mathbb{Z}, z\in \Lambda.
\end{equation}
In the sequel we refer to $P_{x,\omega}(\cdot)$ as the ``quenched"
law. One also defines the ``annealed" law
 on $\Omega\times
\mathbb{Z}^N$ by:
\begin{equation}
  P_x(\cdot)= \int  P_{x,\omega}(\cdot)  \mathbb{P}(d\omega) \mbox{\ \ for\ }x\in \mathbb{Z}.
\end{equation}
In the rest of the paper, we use  $\mathbb{E}$ corresponding to
$\mathbb{P}$, $E_{x,\omega}$ corresponding to $P_{x,\omega}$ and
$E_x$ corresponding to $P_x$ to denote the expectation respectively.
 Define the shift $T$ on $\Omega$ by relation
$(T\omega)_i=\omega_{i+1}$. Let $$
a_i=\frac{\omega_0(-i)+\cdots+\omega_0(-L)}{\omega_0(1)}, 1\le i\le
L,$$
$$A:=A(0)=\left(
                                           \begin{array}{cccc}
                                             a_1 & \cdots & a_{L-1} & a_L \\
                                             1 & \cdots & 0 & 0 \\
                                             \vdots & \ddots & \vdots & \vdots \\
                                             0 & \cdots & 1 & 0\\
                                           \end{array}
                                         \right).
$$
For $k\ge l$ set $A(k,l)=A(k)\cdots A(l),$ $A(k):=T^kA,$ and for
$l\in \mathbb{Z}$ set $$\delta(l,l+1)=1\mbox{ and } \forall k\ge l,
\delta(k,l)=\langle e_1,A(k)\cdots A(l)e_1\rangle.$$ Note that all
$\delta(k,l)$ defined above are strictly positive.

Define the norm of matrix $A$ by
$$\| A\|=<e_1,Ae_1>.$$ We have $\delta(k,l)=\|A_k\cdots A_l\|.$

It is easy to verify that for $n\ge L,$ $A_nA_{n-1}\cdots A_0\gg0,$ where for a matrix $A,$ $A\gg0$ means that all entries of $A$ are strictly positive. Then one follows from Frobenious theory of positive matrices that there exists a number $\lambda_0$ such that $Ax=\lambda_0 x$ for some $x\in \mathbb R^L$ and $|\lambda|<\lambda_0$ for all other eigenvalues of $A.$ Consequently, $A$ is contracting. Next, suppose $V$ is linear subspace of $\mathbb R^L$ with dimension $1\le d\le L.$ Then for any $v\in V,$ $Av\in \mathbb R^{d+1}.$ Therefore the set $\{A(\omega):\omega\in \text{supp}\mu\}$ is strongly irreducible.

For the definition of contracting set and strongly irreducible set, see \cite{bl85}.
Let $$l(A)=\sup\{\log^+\|A\|, \log^-\|A^{-1}\|\}.$$
By the elliptic condition of $\mu,$ we have that for all $u\in \mathbb R,$ $$\mathbb E(e^{ul(A)})<\infty.$$

Then we have the following facts which could be found in \cite{bl85}.

\noindent\textbf{Facts: }{\it (1) The limit \begin{equation}\label{ly}\gamma_L:=\lim_{n\rto}\frac{1}{n}\mathbb E(\log \|A_{n-1}\cdots A_0\|\end{equation} exists;

\noindent (2) For $u\in \mathbb R $ the limit $$F(u):=\lim_{n\rto}\frac{1}{n}\log \mathbb E(\|A_{n-1}\cdots A_0\|^u)$$ exists and the function $F(\cdot)$ is analytic;

\noindent (3) Consequently, $\{\frac{1}{n}\log\|A_{n-1}\cdots A_0\|\}_{n\ge 1}$ satisfies a large deviation principle. Precisely, for any $\epsilon>0$,
$$\lim_{n\rto}\frac{1}{n}\log\mathbb P(\log \|A_{n-1}\cdots A_0\|>\epsilon)=-I(\epsilon),$$
 $$\lim_{n\rto}\frac{1}{n}\log\mathbb P(\log \|A_{n-1}\cdots A_0\|<-\epsilon)=-I(-\epsilon),$$
 where the rate function $I(x)=\sup_{u\in \mathbb R}\{ux-F(u)\}.$}

The number $\gamma_L$ is called  the greatest Liapounov exponent of $A.$ It serves as a criteria for RWRE with bounded jumps. The
following results can be found in Br\'{e}mont \cite{B02},(see page
1271, lemma 4 in page 1272, theorem 2.4 in page 1275, theorem
3.5 in page 1284 respectively). \\
 {\bf Theorem A. } ( Br\'{e}mont,
\cite{B02}) {\it For the RWRE with bounded jumps $X_n$, we have
\begin{itemize}
  \item[1.] There exists a unique unit random vector $V$ with
  strictly positive components and a unique random variable $\lambda$ such that $AV=\lambda T V.$
\item[2.] $\gamma_L=\mathbb{E}(\log\lambda).$
\item[3.]  There exists a constant $C>0$ such that for $k>l$
\begin{equation}\label{es1}
(1/C)(T^k\lambda\cdots T^l\lambda)\le\delta(k,l)\le
C(T^k\lambda\cdots T^l\lambda).
\end{equation}
\item[4.] If $\gamma_L<0,$ then $X_n\rto$ $P_0$-a.s.. If $\gamma_L>0,$ then $X_n\rightarrow-\infty$ $P_0$-a.s.. If $\gamma_L=0,$ the walk is recurrent almost surely.

 \item[5.]If $\mathbb{E}(\sum_{n=1}^\infty T^{n-1}\lambda\cdots \lambda)<\infty,$ then
$\frac{X_n}{n}\rightarrow c>0$ $P_0$-a.s.. If
$\mathbb{E}(\sum_{n=1}^\infty (T^{n-1}\lambda\cdots \lambda)^{-1})<\infty,$ then $\frac{X_n}{n}\rightarrow c<0$
$P_0$-a.s.. If  $\mathbb{E}(\sum_{n=1}^\infty T^{n-1}\lambda\cdots \lambda)=\infty$ and $\mathbb{E}(\sum_{n=1}^\infty (T^{n-1}\lambda\cdots \lambda)^{-1})<\infty,$ then $\frac{X_n}{n}\rightarrow c=0$
$P_0$-a.s.. \qed
\end{itemize}}

We have from (\ref{es1}) that \begin{equation*}
  \lim_{n\rto}\frac{1}{n}\log \mathbb E((T^{n-1}\lambda\cdots\lambda)^u)=\lim_{n\rto}\frac{1}{n}\log \mathbb E(\|A_{n-1}\cdots A_0\|^u)=F(u).
\end{equation*}
Consequently $\mathbb{E}(\sum_{n=1}^\infty T^{n-1}\lambda\cdots \lambda)<\infty$ if and only if $F(1)<0,$ while $\mathbb{E}((\sum_{n=1}^\infty T^{n-1}\lambda\cdots \lambda)^{-1})<\infty$ if and only if $F(-1)<0.$

In this point of view, we have that \begin{align*}
 & F(1)<0\Rightarrow \frac{X}{n}\rightarrow c>0;\\
 & F(-1)<0\Rightarrow \frac{X}{n}\rightarrow c<0;\\
 & F(1)\ge 0, F(-1)\ge0 \Rightarrow \frac{X}{n}\rightarrow c=0.\\
\end{align*}

{\bf Remark 1.1}
\begin{itemize}
  \item[1.] By the convexity of function ${x^{-1}}$,
  we have $F(-1)\ge -F(1).$ As a consequence, there is one and only one case in 5
  of  Theorem A happens.
 \item[2.] Note that the function $\log x$ is concave, we have
 $\gamma_L\le F(1).$ If $\gamma_L
 <0$ ( by 4 of theorem A, $X_n\to\infty$, a.s.) it is possible that
 $F(1)>0.$ Then  by 5 of theorem A we have
$\frac{X_n}{n}\rightarrow 0$  a.s..
Similarly it is also possible that $F(-1)>0$ and $\gamma_L>0.$ In this case, $X_n\to\infty,$ a.s. but $\frac{X_n}{n}\rightarrow 0$ a.s.. In these two situations,  the
slowdown properties occur and we show in the following main theorem that $X_n$ grows with only sub-linear speed.
\end{itemize}

\begin{center}
{\section{\hspace{-0.4cm.}\hspace{0.2cm}Main result and proofs}}
\end{center}
\begin{theorem}For the RWRE with bounded jumps $X_n$,  one of
the following conditions holds

\noindent(1) $\gamma_L:=\lim_{n\rto}\frac{1}{n}\mathbb E(\log\|A_{n-1}\cdots A_0\|)<0$ and $F(1):=\lim_{n\rto}\frac{1}{n}\log \mathbb E(\|A_{n-1}\cdots A_0\|)>0$;

\noindent(2) $\gamma_L:=\lim_{n\rto}\frac{1}{n}\mathbb E(\log\|A_{n-1}\cdots A_0\|)>0$ and $F(-1):=\lim_{n\rto}\frac{1}{n}\log \mathbb E(\|A_{n-1}\cdots A_0\|^{-1})>0.$\\
  Then there exists $s\in(0,1)$ such that
  $\frac{X_n}{n^{s'}}\rightarrow
  0,$  $P_0$-a.s., for all $s'>s.$
\end{theorem}
{\bf Remark 2.1}
\begin{itemize}
  \item[1.] By the discussion in Remark 1.1, $X_n \rto \mbox{ but } \frac{X_n}{n}\rightarrow
  0,$ $P_0$-a.s., while  case (1) of the theorem happens; and similarly  $X_n \to -\infty \mbox{ but }
\frac{X_n}{n}\rightarrow 0, $
  $P_0$-a.s., while case (2) of the theorem happens.
  \item[2.] For a fixed environment $\omega$, we call $[-K\log n, K\log n]$ a ``trap" for the walk if
  the walk with positive (quenched) probability spends more than $n$ steps in $[-K\log n, K\log
  n]$, i.e.,    $P_{0,\omega}\z[T_{[-K\log n, K\log n]}>n\y]\ge
\varepsilon>0$; and we say the fixed environment $\omega$ formulated
a trap for the walk, where $T_{[-K\log n, K\log n]}$ is the first
exit time of the walk from $[-K\log n, K\log n]$.
\item[3.] The key step in the proof of the Theorem is to show that
there are ``many" environments formulated traps for the walk  in the
sense of the following (16) under the conditions of the Theorem. For
this purpose we need some estimations for the (quenched) exit
probabilities and the environment factors.

\end{itemize}

{\noindent\bf Proof of Theorem 2.1:}

Consider integers $(a,b,k)$ with $a<b$ and define
$$P_{k,\omega}\{a,b,-\}:=P_{k,\omega}\{\mbox{the walk reaches }(-\infty, a]\mbox{ before }[b,\infty)\}$$
and similarly $$P_{k,\omega}\{a,b,+\}:=P_{k,\omega}\{\mbox{the walk
reaches }[b,\infty)\mbox{ before }(-\infty, a]\}.$$ We have the
following lemma which was proved in  Br\'{e}mont \cite{B02}.
\begin{lemma}
  If $a<k<b,$ then
  \begin{equation}
P_{k,\omega}\{a,b,-\}=\frac{\sum_{j=k}^{b-1}\delta(j,a+1)}{\sum_{j=a}^{b-1}\delta(j,a+1)}.
  \end{equation}
\end{lemma}

With this Lemma in hands, defining for $z\in\mathbb{Z}$
$$H_z^l:=\min[n>0: X_n\le z],\  H_z^r:=\min[n>0: X_n\ge z],$$
 for $M>1,$ $N>L,$ $1\le k\le L,$  we
have
\begin{eqnarray*}
    &&P_{1,\omega}\left\{H_0^l<H_{M+1}^r\right\}=P_{1,\omega}\{0,M+1,-\}=\frac{\sum_{j=1}^{M}\delta(j,1)}{\sum_{j=0}^{M}\delta(j,1)}=1-\frac{\delta(0,1)}{\sum_{j=0}^{M}\delta(j,1)} \nonumber\\
  &&\quad\quad=1-\frac{1}{\sum_{j=0}^{M}\delta(j,1)}
\ge 1-\frac{1}{\delta(M,1)}=1-e^{-\log \delta(M,1)},
\end{eqnarray*}
and
\begin{eqnarray*}
    &&P_{-k,\omega}\left\{H_0^r<H_{-(N+1)}^l\right\}=P_{-k,\omega}\{-(N+1),0,+\}=1-P_{-k,\omega}\left\{-(N+1),0,-\right\} \nonumber\\
    &&\quad\quad=1-\frac{\sum_{j=-k}^{-1}\delta(j,-N)}{\sum_{j=-(N+1)}^{-1}\delta(j,-N)}
   =1-\frac{\sum_{j=-k}^{-1}\delta(j,-N)}{1+\sum_{j=-N}^{-1}\delta(j,-N)}\nonumber\\
&&\quad\quad\ge 1-\sum_{j=-k}^{-1}\delta(j,-N)\ge
1-\sum_{j=-L}^{-1}\delta(j,-N)\nonumber\\
&&\quad\quad=1-e^{\log \sum_{j=-L}^{-1}\delta(j,-N)}.
\end{eqnarray*}
Let $R_{M}=\frac{1}{M}\log \delta(M,1)$ and let
$R_{N}=\frac{1}{N}\log \sum_{j=-L}^{-1}\delta(j,-N).$ Then we have
\begin{equation}\label{jia}
P_{1,\omega}\left\{H_0^l<H_{M+1}^r\right\}\ge\z(1-e^{-MR_{M}}\y)_+
\end{equation}
and
\begin{equation}\label{jian}
P_{-k,\omega}\left\{H_0^r<H_{-(N+1)}^l\right\}\ge
\z(1-e^{NR_{N}}\y)_+.\end{equation}

From (\ref{ly}) one follows that $\mathbb P$-a.s., $$\lim_{M\rto}R_M=\lim_{N\rto}R_n=\gamma_L.$$

{\noindent \bf Case 1.} Suppose that $\gamma_L<0$ but $F(1)>0.$
Not that $F$ is a strictly convex function satisfying $F(0)=0,$ $F(1)>0$
and $F'(0)<0.$ Therefore, there exists a unique $s\in(0,1)$ such
that $F(s)=0.$ We fix such $s$ in the
remainder of the proof.

We now set, for $U=[-N,M],$
$\gamma(U)=1\wedge\max[e^{NR^-},e^{-MR^+}],$ where $R^-=R_{N}$ and
$R^+=R_{M}.$ Define $T_U:=\inf[k, X_k\in U^c].$ Note that $T_U$ is
the exit time of the walk from the set $U.$ Then we have from the
strong Markov property that
\begin{equation}\label{hit}
    P_{0,\omega}\left(T_U>n\right)\ge P_{0,\omega}\big[\#\{1< k\le T_U,\
    X_{k-1}X_k\le0\}>n\big]\ge\left(1-\gamma(U)\right)^n.
\end{equation}
The first inequality of the last expression follows immediately. For
the second one, we define $$\tilde{H}_0:=\inf[k>1,\
X_kX_{k-1}\le0].$$ Note that $\tilde{H}_0$ can be explained as the
first time the walk crosses 0 after time 1. By decomposing the event
$\{\tilde{H}_0<T_U\}$ according to the value of $X_1,$ we
have\begin{eqnarray}
   &&P_{0,\omega}\left(\tilde{H}_0<T_U\right)=P_{1,\omega}\z(H_0^l<H_{M+1}^r\y)P_{0,\omega}(X_1=1)\nonumber\\
&&\quad\quad+\sum_{j=1}^{L}P_{-j,\omega}\z(H_0^r<H_{-(N+1)}^l\y)P_{0,\omega}(X_1=-j)\nonumber\\
&&\hspace{2cm}\mbox{by the estimation in (\ref{jia}) and (\ref{jian})}\no\\
  && \quad\quad\ge P_{0,\omega}(X_1=1)\z(1-e^{-MR^+}\y)_++\sum_{j=1}^{L}P_{0,\omega}(X_1=-j)\z(1-e^{NR^-}\y)_+\no\\
  &&\quad\quad\ge  P_{0,\omega}(X_1=1)\z(1-\gamma(U)\y)+\sum_{j=1}^{L}P_{0,\omega}(X_1=-j)\z(1-\gamma(U)\y)=1-\gamma(U).\no
\end{eqnarray}
Then (\ref{hit}) follows. The remainder of the proof is similar as
Sznitman \cite{S04}, to make the proof complete we still give the
details here. In particular, if $\gamma_L(U)\le\frac{1}{n},$ we have
\begin{equation*}
P_{0,\omega}\left(T_U>n\right)\ge \left(1-\gamma(U)\right)^n\ge
\z(1-\frac{1}{n}\y)^n\rightarrow e^{-1}.
\end{equation*}
Hence for $n$ large enough, $P_{0,\omega}\left(T_U>n\right)\ge
c=e^{-2}>0.$ Note that for $N\ge
\frac{2}{|\gamma_L|}\log n,$ $M,\ \epsilon$ with $M\epsilon\ge \log n,$ by
independence of $R^+$ and $R^-$ under $\mathbb{P},$
\begin{equation}\label{rr}
    \mathbb{P}\z(\gamma(U)\le\frac{1}{n}\y)\ge  \mathbb{P}\z(R^-\le\frac{\gamma_L}{2},R^+\ge
    \epsilon\y)=\mathbb{P}\z(R^-\le\frac{\gamma_L}{2}\y)\mathbb{P}\z(R^+\ge
    \epsilon\y).
\end{equation}
Then for $n$ large enough and $\eta>0$ small, by the large deviations, we have that
\begin{equation*}
 \mathbb{P}\z(\gamma(U)\le\frac{1}{n}\y)\ge \frac{1}{2}\mathbb{P}\z(R^+\ge
    \epsilon\y)\ge\frac{1}{2}\exp\{-I(\epsilon)M(1+\eta)\}.
\end{equation*}
Now we optimize $\epsilon,$ $M$ by looking at
\begin{equation*}
    \inf[I(\epsilon)M, M\epsilon\ge \log n]= \inf_{\epsilon>0}\z[\frac{I(\epsilon)}{\epsilon}\log
    n\y],
\end{equation*}
and recall that $F(u)=\sup_x[xu-I(x)].$ Let
$\alpha:=\inf_{\epsilon>0}\z[\frac{I(\epsilon)}{\epsilon}\y].$ By a
duality argument(see \cite{DZ98} lemma 4.5.8), we see that
$F(\alpha)=0.$ In the other words
$\inf_{\epsilon>0}\z[\frac{I(\epsilon)}{\epsilon}\y]=s,$ recalling
that $s\in(0,1)$ is the unique positive zero of the function
$F(\cdot).$ Therefore, choosing $K>0,$ $\eta>0$ properly, for large
n, from the discussion above we have
\begin{equation}\label{tr}
\mathbb{P}\z(P_{0,\omega}\z[T_{[-K\log n, K\log n]}>n\y]\ge
e^{-2}\y)\ge n^{-s(1+\eta)}.
\end{equation}
With (\ref{tr}) we have created a trap of size $2K\log n$ which
retains the walk for $n$ units of time with large probability. If
$s'>s,$ choosing $\eta$ small in (\ref{tr}) there will be many such
traps in $[0,n^{s'}]$ which will prevent the walk from moving to
distance $n^{s'}$ from the origin before time $n$. Precisely, for
large $n,$ with $M$ the number of traps in $[0,n^{s'}]$, which is of
order $\frac{n^{s'}}{\log n},$ and with $T_i$ the time to exit the
$i.$th trap after reaching its center, for $\lambda>0,$ we have
\begin{eqnarray*}
 &&P_0\z(X_n>n^{s'}\y)\le P_0\z(T_1+\cdots+T_M<n, X_.\mbox{ reaches the center of the $i.$th trap,
 $i$=1,...,M}\y)\no\\
 &&\quad\quad\le e^{\lambda n}E_0\z(e^{-\lambda(T_1+\cdots+T_M)},X_.\mbox{ reaches the center of the $i.$th trap,
 $i$=1,...,M}\y)\no\\
 &&\quad\quad\quad (\mbox{{\it using Markov property under $P_{0,\omega}$, the independence under $\mathbb{P}$ }}\no\\
 &&\quad\quad\quad \mbox{{\it of the environments in different traps and the
 stationarity}})\no\\
 &&\quad\quad =e^{\lambda n}E_0\z(e^{-\lambda T_{[-K\log n, K\log
 n]}}\y)^M\no\\
 &&\quad\quad \le e^{\lambda n}\z(1-\frac{e^{-2}}{n^{s(1+\eta)}}+\frac{e^{-2-\lambda
 n}}{n^{s(1+\eta)}}\y)^M,\mbox{using the fact $1-x\le
 e^{-x}$}\no\\
 &&\quad\quad\le e^{\lambda n-M\frac{e^{-2}}{n^{s(1+\eta)}}(1-e^{-\lambda
 n})}.
\end{eqnarray*}
Since $M\sim const\frac{n^{s'}}{\log n},$ if we now choose $\eta$
small, $\lambda=\frac{1}{2}n^{s'-s(1+2\eta)-1},$ for $n$ large, we
have
$$P_0\z(X_n>n^{s'}\y)\le e^{\frac{1}{2}n^{s'-s(1+2\eta)}-n^{s'-s(1+2\eta)}}=e^{-\frac{1}{2}n^{s'-s(1+2\eta)}}.$$
Then it follows from Borel-Cantelli lemma that,
 $$P_0\mbox{-a.s., } \lim_{n\rto}\frac{X_n}{n^{s'}}=0, \mbox{ for all $s'>s.$}$$

{\noindent\bf Case 2:} By assumption $\gamma>0,$
$F(-1)>0.$  Then
$F$ is a strictly convex function satisfying $F(0)=0,$ $F(-1)>0$ and
$F'(0)>0.$ Therefore, there exists a unique $s\in(-1,0)$ such that
$F(s)=\log \mathbb{E}\left(\lambda^s\right)=0,$ that is,
$\mathbb{E}\left(\lambda^s\right)=1.$ Fix such $s.$ The proof moves
on as that of  Case 1. Using large deviation and  changing he role of $R^+$ and $R^-$ in (\ref{rr}), we
can get an estimation as (\ref{tr}), i.e.,
\begin{equation}\label{tr1}
\mathbb{P}\z(P_{0,\omega}\z[T_{[-K\log n, K\log n]}>n\y]\ge
e^{-2}\y)\ge n^{s(1+\eta)},
\end{equation} for $n$ large and $\eta>0$ small.
Recall that $s\in  (-1,0)$ in this case. Using a similar argument of
 Case 1  below (\ref{tr}), we can get
$$P_0\mbox{-a.s., } \lim_{n\rto}\frac{X_n}{n^{-s'}}=0, \mbox{ for
all $s'<s,$}$$
 which completes the proof. \qed

\begin{center}
{\section*{References}}
\end{center}

\end{document}